\date{}
\title{Recurrence of Distributional Limits \\ of Finite Planar Graphs}
\author{Itai Benjamini \and Oded Schramm}

\documentclass[12pt,naturalnames]{article}
\usepackage{amsmath}
\usepackage{amsthm}
\usepackage{amsfonts}
\usepackage{graphicx}
\input labelfig.tex
\newif\ifhyper\IfFileExists{hyperref.sty}{\hypertrue}{\hyperfalse}
\ifhyper\usepackage{hyperref}\fi

\newif\ifdraft
\drafttrue
\numberwithin{equation}{section}
\numberwithin{figure}{section}

\newtheorem{theorem}{Theorem}
\numberwithin{theorem}{section}
\newtheorem{corollary}[theorem]{Corollary}
\newtheorem{lemma}[theorem]{Lemma}
\newtheorem{prop}[theorem]{Proposition}

\newtheorem{conjecture}[theorem]{Conjecture}
\newtheorem{problem}[theorem]{Problem}
\def\eref#1{(\ref{#1})}

\newcommand{\R}{\mathbb{R}}

\newcommand{\Z}{\mathbb{Z}}
\newcommand{\N}{\mathbb{N}}

\def\diam{\mathrm{diam}}

\def\ev#1{\mathcal{#1}}

\def \eps {\epsilon}

\def \P {{\bf P}}
\def\md{\mid}
\def\Bb#1#2{{\def\md{\bigm| }#1\bigl[#2\bigr]}}
\def\BB#1#2{{\def\md{\Bigm| }#1\Bigl[#2\Bigr]}}

\def\Pb{\Bb\P}

\def\EB{\BB\E}

\def \E {{\bf E}}

\def \proof {{ \medbreak \noindent {\bf Proof.} }}
\def\proofof#1{{ \medbreak \noindent {\bf Proof of #1.} }}
\def\proofcont#1{{ \medbreak \noindent {\bf Proof of #1, continued.} }}

\def\Schaeffer{MR2001i:05142}

\def\CdVcp{CdV91a}
\def\McCaughan{MR99b:52041}
\def\HStype{MR96h:52017}
\def\RodinSullivan{MR90c:30007}
\def\KoebeCP{pK36}
\def\ADJBooK{MR98i:82001}
\def\Bw95{MR96b:20046}
\def\AAIJW{MR99m:83059}
\def\ABNRW{MR99c:83026}
\def\Ma99{MR2000k:83024}
\def\GW00{MR1741015}
\def\Tu62{MR24:A695}
\def\oHmtp{MR98f:60207}
\def\BLPSgip{MR99m:60149}
\def\BabaiGr{MR1447704}
\def\Thomas{MR2001e:05118}

\def\SS{\mathfrak{S}}

\begin{document}
\maketitle

\begin {abstract}
Suppose that $G_j$ is a sequence of finite connected planar graphs,
and in each $G_j$ a special vertex, called the root, is chosen
randomly-uniformly.
We introduce the notion of a distributional limit $G$ of such graphs.
Assume that the vertex degrees of the vertices in $G_j$ are bounded, and
the bound does not depend on $j$.  Then after passing to a subsequence,
the limit exists, and is a random rooted graph $G$.
We prove that with probability one $G$ is recurrent.
The proof involves the Circle Packing Theorem.
The motivation for this work comes from the theory of random spherical
triangulations.
\end {abstract}

\section{Introduction}

\subsection{Random triangulations}

In recent years, physicists were interested in the study of random
surfaces \cite{\ADJBooK}. Random triangulations turned out to be a useful model
for exact calculations, non rigorous arguments, and Monte-Carlo simulations
regarding the geometry of random surfaces and the behaviour of physical systems
on these surfaces. From a mathematical viewpoint, natural measures that were
considered are the uniform measures on  isomorphism classes of
triangulations of the sphere with a fixed number of vertices. In \cite{\AAIJW,\ABNRW}
diffusion on some random surfaces and random walks on random triangulations including
the uniform measure have been
considered. It was suggested there that the probability
for the random walk to be at time $t$ at its starting vertex should
decay like $t^{-1}$, provided that
$t$ is not too large relative to the size of the triangulation, and that
the mean square displacement at time $t$ is $t^{1/2}$.
Motivated by these observations, we decided to study the recurrence
versus transience dichotomy for
limits of random rooted spherical triangulations.  It turned out
that in the end the results apply to the more general setting of
planar bounded-degree graphs.
It will be proven that under the assumption of a
uniform bound on the vertex degrees, limits of finite
planar graphs are recurrent (provided that the root is chosen uniformly).

\subsection{Limits of graphs}

In addition to studying asymptotic properties of large random objects,
it is mathematically appealing and natural to introduce a limiting infinite
object and study its properties.
In order to define the limit of a sequence of (possibly random) triangulations
or graphs, it is necessary to keep track of a basepoint; or a {\bf root}.

A {\bf rooted graph} is just a pair $(G,o)$, where
$G$ is a (connected) graph and $o$ is a vertex in $G$.
A rooted graph $(G,o)$ is isomorphic to $(G',o')$
if there is an isomorphism of $G$ onto $G'$ which
takes $o$ to $o'$.
In this paper, we only consider {\bf locally finite} graphs;
that is, each vertex has finitely many neighbors.

The space $\mathcal X$ of isomorphism classes of rooted connected
(locally finite) graphs has a natural topology,
which is induced by the following metric.
Let $(G,o),(G',o')\in\mathcal X$. 
For $r=1,2,\dots$, let $B_G(o,r)$ be the closed ball of radius $r$ about $o$
in $G$, and similarly for $G'$.
Let $k$ be the supremum of all $r$ such that $\bigl(B_G(o,r),o\bigr)$ and
$\bigl(B_{G'}(o',r),o'\bigr)$ are
isomorphic as rooted graphs, and set
$d\bigl((G,o),(G',o')\bigr):=2^{-k}$, where $2^{-\infty}:=0$.
Then $d$ is a metric on $\mathcal X$.
Moreover, it is easy to verify that for every
$M<\infty$ the subspace $\mathcal X_M\subset\mathcal X$ of graphs 
with degrees bounded by $M$ is compact in this topology.

Suppose that $H$ is a finite connected graph.  We cannot
think of $H$ as an element of $\mathcal X$, unless a root $o$
is chosen.  The most natural way to choose a root
is to make the choice random, and uniform among the
vertices of $H$.  In this way, a finite unrooted graph
$H$ corresponds to a probability measure $\mu_H$ on
$\mathcal X$.  More explicitly, for every Borel subset
$\ev A\subset\mathcal X$, $\mu_H(\ev A)$ is
equal to the probability that $(H,o)\in\ev A$ when $o$
is chosen uniformly and randomly among the vertices of $H$.

Suppose that $(G,o)$ is a random rooted finite graph.
(This means that $(G,o)$ is a sample from a
Borel probability measure $\mu$ on $\mathcal X$,
which is called the {\bf law} of $(G,o)$,
and $\mu$ is supported on the set of finite graphs.)
Then $(G,o)$ is {\bf unbiased} if its law is in the
closed convex hull of the measures $\mu_H$.
In other words, for every finite graph $H$,
conditioned on the event that $G$ is isomorphic to
$H$, the distribution of $(G,o)$ is $\mu_H$
(provided that $G$ is isomorphic to $H$ with positive
probability).
Informally, a random rooted finite graph $(G,o)$ is unbiased,
if given $G$ the root $o$ is uniformly distributed among
the vertices $V(G)$.

Let $(G,o)$ and $(G_1,o_1),(G_2,o_2),\dots$
be random connected rooted graphs.
We say that $(G,o)$ is the {\bf distributional limit} of $(G_j,o_j)$
as $j\to\infty$ if for every $r>0$ and for every finite rooted graph $(H,o')$,
the probability that $(H,o')$ is isomorphic to
$\bigl(B_{G_j}(o_j,r),o_j\bigr)$
converges to the probability that
$(H,o')$ is isomorphic to $\bigl(B_{G}(o,r),o\bigr)$.
This is equivalent to saying that the law of $(G,o)$,
which is a probability measure on $\mathcal X$, is the weak limit
of the law of $(G_j,o_j)$ as $j\to\infty$.

It is easy to see, by compactness or diagonalization, 
that if $(G_j,o_j)\in\mathcal X_M$ a.s., $M<\infty$, then
there is always a subsequence of the sequence $(G_j,o_j)$
having a distributional limit.

\begin{theorem}\label{para}
Let $M<\infty$, and
let $(G,o)$ be a distributional limit
of rooted random unbiased finite planar graphs $G_j$
with degrees bounded by $M$.
Then with probability one $G$ is recurrent.
\end{theorem}

\medskip
\noindent{\bf Remarks}.
To illustrate the theorem, the reader may wish to consider the case where
each $G_j$ is a finite binary tree of depth $j$.  (In that case
the distance of the root from the leaves of the tree is approximately
a geometric random variable.  Hence, $G$ is not the $3$-regular tree.)
\par
The assumption of planarity in the theorem is necessary, as can be seen by
considering the intersection of $\Z^3$ with larger and larger balls.
In this case, since the surface area to volume ratio tends to zero,
the root is not likely to be close to the boundary.
The distributional limit will then be $(\Z^3,0)$ a.s.
Also, the $3$-regular tree can be obtained as
the a.s.\ limit of unbiased finite graphs ($3$ regular
graphs with girth going to infinity).
\par
A natural extension of the collection of all planar
graphs is the collection of all graphs with an excluded minor.
It is reasonable to guess that if $H$ is any graph
and one replaces the assumption of planarity by the assumption
that all $G_j$ do not have $H$ as a minor, then the theorem still
holds, because many theorems on planar graphs generalize to excluded
minor graphs (see, e.g.,~\cite{\Thomas}).
\par
Because of the assumption of a uniform bound on the degrees,
the interesting case when $G_n$ is uniformly distributed
among all isomorphism classes of triangulations of size $n$ is
unfortunately not included.
We {\bf conjecture} that the theorem still holds for limits of these measures.
(See \cite{\Ma99,\GW00}
for a study of the degree distribution for those measures.)
If $G_n$ is uniformly distributed among spherical triangulations
of size $n$ with degree at most $M$, then the theorem above does
apply to any (subsequential) limit of $G_n$.
\par 
A key to almost all the rigorous results regarding random
surfaces and triangulations
is the enumeration techniques, which originated with the fundamental work
of Tutte \cite{\Tu62}.
Schaeffer~\cite{\Schaeffer} found a simpler proof for Tutte's enumeration,
and thereby produced a good sampling algorithm.
  In contrast, our approach in this paper makes no use of enumeration.
The proof of \ref{para} is based on the theory of circle packings with specified
combinatorics.
\medskip 

Following is another statement of the theorem from a slightly
different perspective.  Let $G$ be a finite graph
and let $X_0,X_1,\dots$ be simple random walk
on $G$ started from  a random-uniform vertex $X_0\in V(G)$.
Let $\phi(n,G)$ be the probability that $X_j\ne X_0$
for all $j=1,2,\dots,n$, and set
$$
\phi_M(n):=\sup \bigl\{\phi(n,G) :G\in \mathcal{G}_M\bigr\},
$$
where $\mathcal{G}_M$ denotes the collection of finite
planar graphs with maximum degree at most $M$.

\begin{corollary} For all $M<\infty$,
$$
\lim_{n\to\infty} \phi_M(n)=0.
$$
\end{corollary}
\proof
Let $h:=\inf_{n\in\N} \phi_M(n)$.
Then there is a sequence
of graphs $G_j\in\mathcal{G}_M$ with $\phi(n,G_j)\ge \phi(j,G_j)\ge h/2$
for all $n\le j$.
Let $(G,o)$ be a distributional limit of a subsequence of
$(G_j,o_j)$, where $o_j$ is unbiased in $V(G_j)$.
Then for every $n\in\N$ the probability that simple
random walk on $G$
started from $o$ will not revisit $o$ in the first
$n$ steps is at least $h/2$.  By the theorem, $G$ is
recurrent, and hence $h=0$, as required.
\qed

\begin{problem}
Determine the rate of decay of $\phi_M(n)$ as a function of $n$.
\end{problem}

The example of an $n\times n$ square in the
grid $\Z^2$ shows that $\inf_n\phi_M(n)\log n>0$.
It might be reasonable to guess that $\phi_M(n)$ decays like $c(M)/\log n$.

\medskip

In the next section, a proof of Theorem~\ref{para} is given,
and the last section is devoted to miscellaneous remarks,
including examples of planar triangulations with
uniform growth $r^\alpha$, $\alpha\notin \{1,2\}$.

\section{Recurrence}

Theorem \ref{para} will be proved using the theory of combinatorially specified
circle packings.  At the foundation of this theory is
the Circle-Packing Theorem, which states
that for every finite planar graph $G$ there is
a disk-packing $P$ in the plane whose tangency graph is $G$;
which means that the disks in $P$ are indexed by the vertices
$V(G)$ of $G$ and two disks are tangent iff
the corresponding vertices share an edge.
When $G$ is the 1-skeleton of a triangulation of
the sphere, the packing $P$ is unique up to
M\"obius transformations.
The Circle Packing Theorem was first proved by Koebe
\cite{\KoebeCP} and later various generalizations
have been obtained.  See, e.g., \cite{\CdVcp}
for a proof of this result.
The Circle Packing Theorem
has been instrumental in answering some problems
about the geometry and potential theory of
planar graphs \cite{MP92b,MT90,BSc96,JS00}.
(The relations between circle-packings and analytic function theory
are important, fascinating, and much studied, but are
not as relevant to the present paper.)
Some surveys of the theory of combinatorially specified
circle packings are also available \cite{mB92,hSa94,kS99}. 

The major step in the proof of Theorem~\ref{para} is
the case of triangulations, that is,

\begin{prop}[The triangulation case]\label{triang}
Let $M<\infty$, and
let $(T,o)$ be a distributional limit
of rooted random unbiased (finite) triangulations of the sphere $T_j$
with degrees bounded by $M$.
Then with probability one $T$ is recurrent.
\end{prop}

\proof
We assume, as we may, that $T$ is a.s.\ infinite.
By the Circle-Packing Theorem,
for each $j$, there is a disk-packing $P^j$ in the plane
with tangency graph $T_j$.  Since $T^j$ is random,
also $P^j$ is random. (Actually, the randomness in $T^j$ plays no
role in the proof, and we could assume that $T^j$ is deterministic.
The essential randomness is that of $o_j$.)
We choose $P^j$ to be independent of $o_j$ given $T^j$.
Our immediate goal is to take an appropriate limit
of the disk-packings $P^j$ to obtain a disk-packing $P$
with tangency graph $T$.

\begin{figure}[h]
\centerline{
\includegraphics*[height=2.5in]{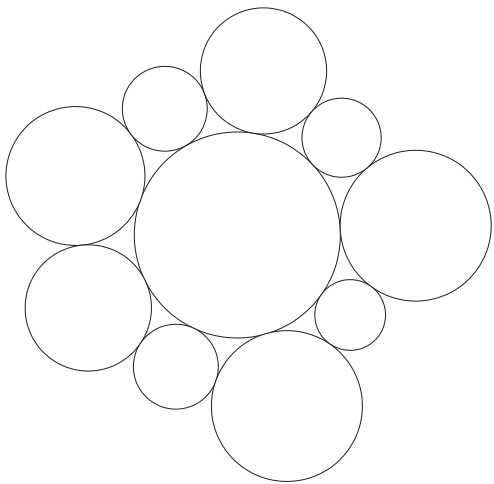}
}
\caption{\label{petals}The Ring Lemma setup.}
\end{figure}

There is a unique triangle $t_j$ in $T_j$
whose vertices correspond to three
disks of $P^j$ which intersect the boundary of the
unbounded component of $\R^2\setminus P^j$.
For every vertex $v\in V(T_j)\setminus t_j$,
the disks corresponding to neighbors of $v$ in $T_j$
surround the disk $P^j_v$ like petals of a flower,
as in Figure~\ref{petals}.
By shrinking the packing, if necessary, assume that
$P^j$ is contained in the unit disk $B(0,1)$.

Let $\hat P^j$ be the image of the packing $P^j$ under the
map $z\mapsto a z+b$, where $a\in(0,\infty)$ and $b\in\R^2$
are chosen so that $\hat P^j_{o_j}$ is the unit disk
$B(0,1)$.

Here is a simple but important fact about disk
packings, known as the Ring Lemma \cite{\RodinSullivan}.
If a disk $P_0$ is surrounded by $n$ other disks
$P_1,\dots,P_n$, as in Figure~\ref{petals},
then the ratio $r_0/r_1$ between the radius of $P_0$
and the radius of $P_1$ is bounded from above by a constant which
depends only on $n$.  Since the vertex degrees in
the triangulations $T_j$ are all bounded by $M$,
it follows that for every $d$ there is some upper bound
$c=c(d,M)$ for the ratio $r/r'$ between any two
radii of disks corresponding to vertices at combinatorial distance
at most  $d$ from $o_j$, provided that $o_j$ is at combinatorial distance
at least $d+1$ from $t_j$.
Because $|V(T_j)|\to\infty$ as $j\to\infty$,
with probability tending to $1$, $P^j_{o_j}$
is not one of the boundary disks corresponding to
vertices of $t_j$.  Moreover, by the uniform bound on the degrees,
the combinatorial
distance in $T^j$ from $o_j$ to $t_j$ tends to
infinity in probability as $j\to\infty$,
because for every $r\ge 1$ the number of vertices of $T^j$ at
distance at most $r$ from $t_j$ is bounded, and hence $o_j$ is
unlikely to be there.
We may therefore conclude that for every $d$
there is a constant $c=c(d)$
so that with probability tending to $1$ all the
disks in $\hat P^j$ corresponding to vertices
at combinatorial distance at most $d$ from $o_j$
have radii in $[1/c,c]$.  By compactness,
we may take a subsequence of $j$ tending to $\infty$
so that the law of $\hat P^j$ tends (in the appropriate
topology) to a random disk packing $P$ whose tangency
graph is $T$.  Assume, with no loss of generality,
that there is no need to pass to a subsequence.

An {\bf accumulation point} of $P$ is a point $z\in\R^2$
such that every open set containing $z$ intersects infinitely
many disks in $P$.  Below, we prove the following result.

\begin{prop}\label{accum}
With probability $1$, there is at most one accumulation point
in $\R^2$ of the packing $P$.
\end{prop}

\proofcont{Proposition \ref{triang}}
  Assuming \ref{accum}, the proof is completed as follows.
In~\cite[Thms.~2.6, 3.1.(1), 8.1]{\HStype}
and independently in~\cite{\McCaughan} it was shown that
if $G$ is a bounded degree tangency graph of a disk-packing
in the plane which has no accumulation points in the plane,
then $G$ is recurrent.
This completes the proof if $P$ has no accumulation points in $\R^2$.
If $P$ has one accumulation point $p\in\R^2$,
then consider the subgraph $G_1$ of $T$ spanned by the vertices
corresponding to disks contained in the disk $B(p,1)$ of radius $1$
about $p$.  By inverting in the circle of radius $1$ about $p$,
it follows that $G_1$ is recurrent.  (Note
that the above quoted results do not require the graph to be
the $1$-skeleton of a triangulation.) 
Similarly, the subgraph $G_2$
spanned by the vertices of $T$ corresponding to disks that
intersect the complement of $B(p,1)$ is also recurrent.
Since $V(T)=V(G_1)\cup V(G_2)$ and the boundary
separating $G_1$ and $G_2$ is finite,
it follows that $T$ is recurrent.
%
\qed
\bigskip

For the proof of Proposition~\ref{accum}, a lemma will be
needed, but some notations must be introduced first.
Suppose that $C\subset\R^2$ is a finite set of points.
(In the application below, $C$ will be the set of centers
of disks in $P^j$.)
When $w\in C$, we define its {\bf isolation radius} as
$\rho_w:=\inf\bigl\{|v-w|:v\in C\setminus\{w\}\bigr\}$.
Given $\delta\in(0,1),\,s>0$ and $w\in C$, we say that $w$
is $(\delta,s)$-{\bf supported} if
in the disk of radius $\delta^{-1}\rho_w$, there are more than
$s$ points of $C$ outside of every disk
of radius $\delta\rho_w$; that is, if
$$
\inf_{ p\in\R^2}\,\,
\Bigl|C\cap B(w,\delta^{-1} \rho_w)\setminus B(p,\delta\rho_w)
\Bigr|\ge s\,.
$$

\begin{lemma}\label{supported}
For every $\delta\in(0,1)$ there is a constant
$c=c(\delta)$ such that
for every finite $C\subset\R^2$ and every $s \ge 2$ the set
of $(\delta,s)$-supported points in $C$
has cardinality at most $c |C|/s$.
\end{lemma}

\proof Let $k\in\{3,4,5,\dots\}$.
Consider a bi-infinite sequence $\SS=(\SS_n:n\in\Z)$ of
square-tilings $\SS_n$ of the plane,
where for all $n\in\Z$ all the squares in
the tiling $\SS_{n+1}$ have the same size and each square in the tiling $\SS_{n+1}$
is tiled by $k^2$ squares in the tiling $\SS_{n}$.
Let $\hat \SS$ denote the collection of all the squares
in all the tilings $\SS_n$.
Assume that no point of $C$ lies on the boundary of a square
in $\hat \SS$.
For every $n\in\Z$,
say that a square $S\in\SS_n$ is $s$-supported
if for every square $S'\in \SS_{n-1}$, we have
$\bigl|C\cap S\setminus S'\bigr|\ge s$.
To estimate the number of $s$-supported squares in $\hat\SS$,
we now define a ``flow'' $f$ on $\hat\SS$.
Set
$$
f(S',S):=\begin{cases}
\min\bigl\{s/2,|S'\cap C|\bigr\}& S\in\SS_{n+1},\ S'\in\SS_n,\ S'\subset S,\\
0 & S\in\SS_{n+1},\ S'\in\SS_n,\ S'\not\subset S,\\
-f(S,S') & S\in\SS_{n},\ S'\in\SS_{n+1},\\
0 & S\in\SS_{n},\ S'\in\SS_m,\ |m-n|\ne 1.\\
\end{cases}
$$
Let $a\in\Z$ be small enough so that each square of
$\SS_a$ contains at most one point of $C$,
and let $b\in\Z,\, b>a$.
Then
\begin{equation*}
\sum_{S'\in\SS_a}\sum_{S\in\SS_{a+1}}f(S',S)=|C|\,,
\qquad
\sum_{S'\in\SS_b}\sum_{S\in\SS_{b+1}}f(S',S)\ge 0\,.
\end{equation*}
Therefore
\begin{equation}\label{teles}
\sum_{n=a+1}^b \,
\sum_{S\in \SS_n}\,
\sum_{S'\in\hat \SS} f(S',S)
\le |C|\,,
\end{equation}
because for every pair $S'\in \SS_n$, $S\in\SS_{n+1}$
with $n\in\{a+1,\dots,b-1\}$, the corresponding two
terms $f(S,S')$ and $f(S',S)$ both appear on
the left hand side and they cancel each other.
By the definition of $f$, for every $S\in\hat\SS$
we have $\sum_{S'\in\hat\SS} f(S',S)\ge 0$.
On the other hand, if $S\in\hat\SS$ is $s$-supported,
then $\sum_{S'\in\hat\SS} f(S',S)\ge s/2$.
Therefore, \eref{teles} implies that the number
of $s$-supported squares in $\hat\SS$ is at most $2|C|/s$.

To estimate the number of $(\delta,s)$-supported
points of $C$, we will compare it with the expected
number of $s$-supported squares in $\hat\SS$,
where $\SS$ is chosen randomly, as follows.
Take $k:=\lceil 20 \delta^{-2}\rceil$ as the parameter
for the sequence $\SS$.
Let $\SS$ have the
distribution such that the distribution of $\hat\SS$ is
invariant under translation and rescaling, and
such that the diameter of any square in $\SS_0$
is in the range $[1,k)$.

To be explicit, we now construct this distribution.
Let $(\alpha_n,n\in\Z)$ be a sequence
of independent random variables with each $\alpha_n$ uniform
in $\{0,1,\dots,k-1\}^2$.  Let $\beta$
be uniform in $[0,\log k)$ and independent from
the sequence $(\alpha_n)$.  Then we may take
$$
\SS_n=
\Bigl\{e^\beta k^n
\bigl([j,j+1]\times[j',j'+1]\bigr)+e^\beta\sum_{m=-\infty}^{n-1}k^{m}\alpha_m
: j,j'\in\Z\Bigr\}.
$$

Let $N$ be the number of $(\delta,s)$-supported points
in $C$.
Say that a point $w\in C$ is a {\bf city} in a square $S\in\hat\SS$
if the edge-length of $S$ is in the range
$[4\delta^{-1}\rho_w,5\delta^{-1}\rho_w]$
and the distance from $w$ to the center of $S$ is at most
$\delta^{-1}\rho_w$.
It is easy to see that there is a constant $c_0=c_0(k)>0$,
which does not depend on $C$ or $w$, such
that $w$ is a city for some square of $\hat\SS$ with probability
at least $c_0$.
By the above choice of $k$, if $w$ is a city in $S$
and $w$ is $(\delta,s)$-supported,
then $S$ is $s$-supported.
Consequently, the expected number of pairs
$(w,S)$ such that $w$ is a city in $S$ and
$S$ is $s$-supported is at least
$c_0 N$.  However, by area considerations it
is clear that there is a constant
$c_1=c_1(\delta)$ such that any square $S$ cannot have
more than $c_1$ cities in it.  Hence,
the expected number of $s$-supported squares
is at least $c_0 N/c_1$.  However, we have seen
above that the number of $s$-supported squares
is at most $2|C|/s$.  Hence $N\le 2 {c_0}^{-1} c_1 |C|/s$,
which proves the lemma.
\qed
\bigskip

\proofof{Proposition \ref{accum}}
Suppose that there is a positive probability that $P$ has
two distinct accumulation points in $\R^2$.
Then there is a $\delta\in (0,1)$ and an $\eps>0$ such that with
probability at least $\eps$ there are two accumulation points
$p_1,p_2$ in $B(0,\delta^{-1})$ such that $|p_1-p_2|\ge 3\delta$.
But this implies that for arbitrarily large $s$ and for infinitely
many $j$ there is probability at least $\eps$ that the center
of $P^j_{o_j}$ is $(\delta,s)$-supported in the set $C^j$ of
centers of the disks in $P^j$, contradicting Lemma~\ref{supported}.
\qed\medskip

\proofof{Theorem \ref{para}}
We claim that there is a constant $c$ such that
for all $j=1,2,\dots$ there is a triangulation $T_j$ of the sphere
with maximum vertex degree at most $c M$, which contains
a subgraph isomorphic to $G_j$, and such that
$|V(T_j)|\le c |V(G_j)|$.  Indeed, embed $G_j$ in the plane,
and let $f$ be a face of this embedding.
Let $v_0,\dots,v_{k-1}$ be the vertices on the boundary
of $f$, in clockwise order.  First, suppose that
$v_j\ne v_{i}$ for $j\ne i$, and that
$v_j$ does not neighbor with $v_i$ unless $|i-j|=1$ or
$\{i,j\}=\{0,k-1\}$.
In that case, we triangulate $f$ by a zigzag pattern;
that is, put into $T_j$ the edges $[v_j,v_{k-j}]$,
$j=1,2,\dots,\lfloor k/2\rfloor-1$
and the edges $[v_j,v_{k-1-j}]$,
$j=1,2,\dots,\lfloor (k-1)/2\rfloor-1$.
Otherwise, first put inside $f$ a cycle of length $k$,
$u_0,\dots,u_{k-1}$, then put the edges
$[u_j,v_j]$, $j=0,\dots,k-1$ and $[u_j,v_{j+1}]$, $j=0,\dots,k-2$
and $[u_{k-1},v_0]$, and then triangulate the face bounded
by the new cycle.
It is easy to verify that when this construction
is applied to every face $f$ of the embedding of $G$, the resulting
triangulation satisfies the required conditions for some appropriate
constant $c$.

Let $o_j$ be a vertex chosen uniformly in  $V(T_j)$.
Since,
$|V(T_j)|\le c |V(G_j)|$, we have $\Pb{o_j\in V(G_j)}\ge 1/c$.
Proposition~\ref{triang} implies that a subsequential limit
of $(T_j,o_j)$ is recurrent a.s.  By Rayleigh monotonicity,
$G$ is recurrent a.s.
\qed

\section{Concluding Remarks}

\subsection{Limits of uniform spherical triangulations}

Let $M\in[6,\infty]$; and let $T^M_j$ be chosen randomly-uniformly among
isomorphism classes of spherical triangulations with $j$ vertices
and maximum degree at most $M$,
and given $T^M_j$ let $o_j$ be chosen uniformly among the vertices
of $T^M_j$.

\begin{conjecture} The distributional limit of $(T^M_j,o_j)$ exists.
\end{conjecture}

This holds when $M=6$, and then the limit is the hexagonal grid
(since by Euler's formula there can be in this case at most $12$
vertices with degree smaller than $6$).
After a first draft of the current paper has been distributed, a proof
of the conjecture for the case $M=\infty$ (that is, no restriction
on the degrees) has been obtained by Omer Angel and Oded Schramm.
A paper with this result is forthcoming.

Assuming the conjecture for now, let $(T^M,o)$ denote the limit random rooted
triangulation.  Let $p_n(T^M)$ denote the probability that the simple random
walk starting from $o$ will be at $o$ at time $n$, given $T^M$.  Following
the discussion in the introduction, it should be believed that $p_n(T^M)$
decays like $n^{-1}$ for almost all $T^M$.
Theorem~\ref{para} shows that when $M<\infty$
the decay cannot be faster than $n^{\alpha}$ with $\alpha<-1$.

\subsection{Intrinsic Mass Transport Principle}

Proposition~\ref{accum} implies that every distributional limit
of unbiased random bounded degree triangulations
of the sphere has at most two ends.
(If $m\in\N$, the statement that a graph $G$ has $m$ ends
is equivalent to the statement that $m$ is the maximum
number of infinite components of $G\setminus W$ as
$W$ ranges over finite subsets of $V(G)$.
The definition of the space of ends is a bit more complicated, and
can be found in many textbooks on point-set topology.)
It is easy to verify that the proof applies also
to the case where there is no uniform bound on the degrees, assuming
that the limit graph is locally finite a.s. 
In fact, one can show that under the assumptions of Theorem~\ref{para},
$G$ has at most two ends.  Moreover, if $G$ is
the distributional limit of unbiased finite rooted graphs,
then a.s.~if $G$ is recurrent it has at most two ends. 
This can be proved using an intrinsic
version of the Mass Transport Principle (MTP).
The extrinsic version of MTP \cite{\oHmtp,\BLPSgip} in its simplest
form says that $\sum_{x\in \Gamma} f(x,y)=\sum_{x\in \Gamma} f(y,x)$ where
$\Gamma$ is a discrete countable group and
$f:\Gamma\times \Gamma\to[0,\infty)$ is invariant under
the diagonal action of $\Gamma$ on $\Gamma\times\Gamma$,
$\gamma:(g,h)\to (\gamma g,\gamma h)$.

We now briefly explain the intrinsic MTP.
Suppose that $f(G,x,y)$ is a function which takes
a graph $G$ and two vertices $x,y$ in $G$ and returns a
non-negative number.  We assume that $f$ is isomorphism
invariant, in the sense that
$f(\psi(G),\psi(x),\psi(y))= f(G,x,y)$,
when $\psi$ is a graph isomorphism.  Also assume
that $f$ is measurable on the appropriate space of
connected graphs with two distinguished vertices.
Then $f$ will be called a {\bf transport function}.
(For example, if $g(G,v,u)$ is the degree of $v$
if $v$ neighbors with $u$ in $G$ and zero otherwise, then $g$ is 
a transport function.)
A probability measure $\mu$ on $\mathcal X$ is said
to satisfy the intrinsic MTP (IMTP) if for every transport function $f$
we have
\begin{equation}\label{e.imtp}
\EB{\sum_{v} f(G,o,v)}=
\EB{\sum_{v} f(G,v,o)}\,,
\end{equation}
where the sum on both sides is over all the vertices of $G$.
It is easy to verify that every unbiased probability measure
on finite graphs in $\mathcal X$ satisfies the IMTP. 
We claim that this also extends to weak limits of such measures. 
Indeed, let $\mu$ be a weak limit of Borel probability measures
$\mu_1,\mu_2,\dots$ on $\mathcal X$ satisfying the IMTP.
Given a transport function $f$ and some $k>0$, let
$f_k(G,v,u):=f(G,v,u)$ when $f\le  k$ and the distance
from $v$ to $u$ in $G$ is at most $k$, and $f_k(G,v,u):=0$, otherwise. 
It is clear that both sides of~\eref{e.imtp} for $\mu$ and $f_k$
are the limits of the corresponding terms for $\mu_j$,
hence they are equal.  The fact that $\mu$ satisfies
the IMTP now follows by taking the supremum with respect to $k$.

Many of the applications of MTP, such as appearing in \cite{\BLPSgip},
for example, can therefore be applied to (limits of finite)
unbiased random graphs, we hope to further pursue this in a future work.

\subsection{Examples of triangulations with uniform growth}

\begin{figure}
\SetLabels
\L(.19*.02)$v_0$\\
\L(.19*.98)$v_1$\\
\endSetLabels
\centerline{\AffixLabels{
\includegraphics*[height=2.5in]{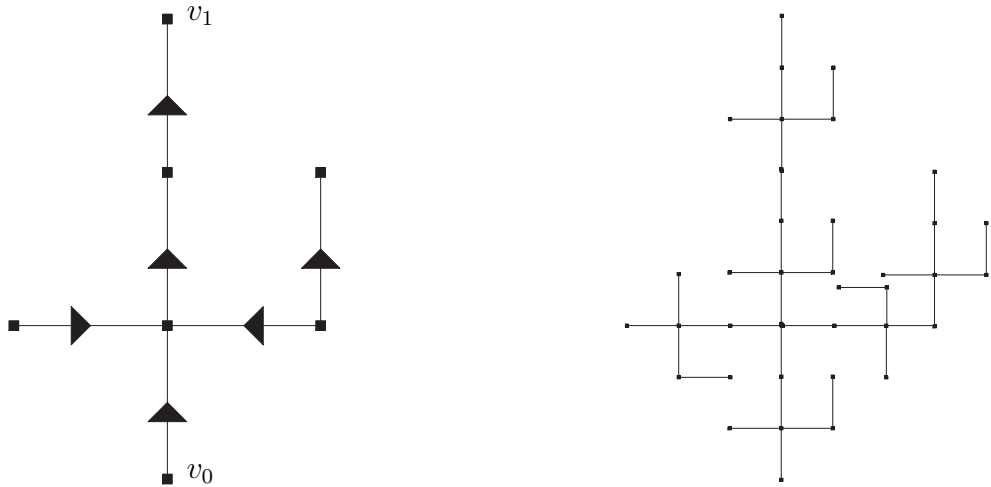}
}}
\caption{\label{ftree}A subdivision rule for a tree.}
\end{figure}

It was observed by physicists \cite{\ADJBooK} that for some random
triangulations
the average volume of balls of radius $r$ is $r^4$. This looks
surprising and calls for an intuitive explanation.
We can't provide one, but to shed  some light on the
geometry of random surfaces, we can construct for every $\alpha>1$ a
triangulation of the plane for which every ball of radius $r$ has $r^\alpha$
vertices, up to a multiplicative constant.

\begin{figure}
\centerline{
\hfil
\includegraphics*[height=2.5in]{ss.eps}
\hfil
\includegraphics*[height=2.5in]{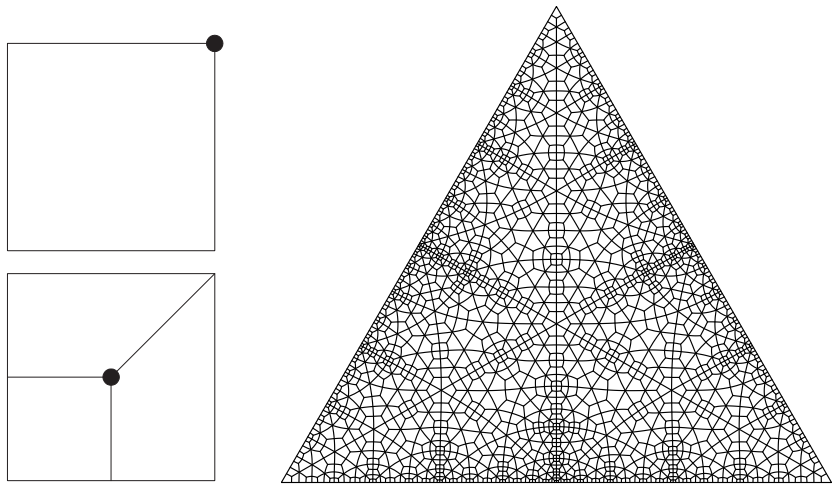}
\hfil
}
\caption{\label{sspt}A subdivision rule and a resulting self-similar tiling.}
\end{figure}

The simplest such construction is based on a tree.
Consider a finite tree $t_1$ with two distinct marked vertices
$v_0$, $v_1$, both having degree 1.  Direct the
edges of $t_1$ arbitrarily.
Let $t_2$ be the tree obtained from $t_1$ by
replacing each directed edge $[u_0,u_1]$ by
a new copy of $t_1$, where $v_0$ replaces $u_0$
and $v_1$ replaces $u_1$.  See Figure~\ref{ftree}
for an example.  Inductively, let $t_n$ be obtained
from $t_{n-1}$ by replacing each edge of $t_{n-1}$
with a copy of $t_1$.
Note that the maximum degree in $t_n$ is the maximum
degree in $t_1$.
Suppose that $t_1$ has $k$ edges and the distance from $v_0$
to $v_1$ in $t_1$ is $\Delta$.  It is then clear that
the diameter of $t_n$ is $\Delta^n$, up to an additive
constant, and the number of edges of $t_n$ is exactly $k^n$.
It follows that for $m \le n$ every ball of radius $\Delta^m$ in $t_n$
 has $k^m$ edges, up to a multiplicative constant, because
$t_n$ is obtained by appropriately replacing each edge of $t_{n-m}$
by a copy of $t_m$.
Then every ball of radius $r\le\diam(t_n)$ in $t_n$
has about $r^{\log k/\log\Delta}$ vertices.
We may pick a root in each $t_n$ and take a subsequential
limit, to obtain a tree $t_\infty$ where every ball
of radius $r$ has about $r^\alpha$ vertices, where $\alpha=\log k/\log\Delta$.

One can easily modify the construction to obtain a similar tree
with growth $r^\alpha$ where $\alpha>1$ is not the ratio of logs of
integers, by letting the replacement rule from $t_{n-1}$ to $t_n$
appropriately depend on $n$.

It is easy to make planar triangulations with similar properties.
For example, suppose that the maximum degree in $t_\infty$
is $M$, and let $T$ be a triangulation of the sphere with at least $M$
disjoint triangles.
Then we may replace each vertex $v$ of
$t_\infty$ by a copy $T_v$ of $T$ and for each edge $[v,u]$
in $t_\infty$ glue a triangle in $T_v$
to a triangle in $T_u$, in such a way that every vertex of
$T_v$ is glued to at most one other vertex.  The resulting
graph is the $1$-skeleton of a planar triangulation, as required.

\medskip

Another example which is somewhat similar but not tree-like appears in
Figure~\ref{sspt}.  It is obtained by starting
with a quadrilateral with a marked corner,
subdividing it as in the figure to obtain
three quadrilaterals with the interior
vertex as the marked corner of each, and continuing inductively.
The result is a map of the plane with quadrilateral faces and
maximum degree $6$.

These examples answer Problem 1.1 from \cite{\BabaiGr}

{\bigskip\noindent\bf Acknowledgment}.
We thank Bertrand Duplantier for introducing us to the
fascinating subject of random spherical triangulations.

\bibliographystyle{halpha}
\bibliography{mr,prep,notmr,kensCP}

\end{document}